\documentclass{amsart}

\usepackage{authblk}
\usepackage{blindtext}
\usepackage{amssymb,mathrsfs,amsmath,dsfont,array}
\usepackage{multirow,arydshln}
\usepackage{tikz}
\usetikzlibrary{arrows,positioning}
\usepackage{capt-of}
\usepackage{mathbbol}     
\usepackage{color}
\usepackage{verbatim} 
\usepackage[all,cmtip]{xy}
\usepackage{fancyhdr}
\usepackage{amsbsy}
\usepackage{amsmath}
\usepackage{nccmath}

\usepackage{hyperref}
\usepackage[multiple]{footmisc}
\usepackage{perpage}
\MakePerPage{footnote}
\renewcommand{\thefootnote}{\fnsymbol{footnote}}

\newtheorem{Thm}{Theorem}[section]
\newtheorem{Pro}[Thm]{Proposition}
\newtheorem{cor}[Thm]{Corollary}
\newtheorem{lem}[Thm]{Lemma}

\newtheorem{deft}[Thm]{Definition}
\newtheorem{rem}[Thm]{Remark}

\numberwithin{equation}{section}

\renewcommand{\tilde}{\widetilde}

\def\a1s{a_1,\cdots, a_s}
\def\a{\alpha}

\def\aa{\mathcal A}

\def\andd{\quad\hbox{and}\quad}

\def\b{\beta}

\def\bl4{B_{\ell\geq4}}

\def\bbbc{{\mathbb C}}

\def\d{\delta}
\def\D{\Delta}

\def\fg{\mathfrak{g}}

\def\hh{{\mathcal H}}

\def\fh{\mathfrak{h}}

\def\lam{\lambda}
\def\Lam{\Lambda}
\def\LL{\mathcal{L}}

\def\ep{\epsilon}
\def\fm{(\cdot,\cdot)}

\def\bbbr{{\mathbb R}}

\def\supp{\hbox{\rm supp}}
\def\sspan{\hbox{\rm span}}

\def\1k{\frac{1}{k}}
\def\op{\oplus}
\def\ot{\otimes}

\def\sub{\subseteq}
\def\sg{\sigma}

\def\pf{\noindent{\bf Proof. }}

\def\bbbz{{\mathbb Z}}

\def\1il{1\leq i\leq\ell}

\begin{document}

\title{Tight irreducible finite weight modules over twisted affine Lie superalgebras}
\thanks{2010 Mathematics Subject Classification: 17B10, 17B65.}
\thanks{Key Words: Twisted Affine Lie Superalgebras; Tight Modules; Finite Weight Modules.}
\bigskip
\chead{Tight irreducible finite weight modules}
\lhead{\thepage}
\rhead{}
\cfoot{}
\maketitle

\centerline{Malihe Yousofzadeh\footnote{Department of Pure Mathematics, Faculty of Mathematics and Statistics, University of Isfahan, Isfahan, Iran,
P.O.Box 81746-73441 and School of Mathematics, Institute for Research in
Fundamental
Sciences (IPM), P.O. Box: 19395-5746, Tehran, Iran. This research  was in part supported by
a grant from IPM (No. 99170216) and  is partially carried out in
IPM-Isfahan Branch.
}}

{\let\thefootnote\relax\footnote{email address: ma.yousofzadeh@sci.ui.ac.ir \& ma.yousofzadeh@ipm.ir.}}

\medskip
\begin{abstract}
For a twisted affine Lie superalgebra with nonzero odd part, we study {tight irreducible  weight modules}   with bounded weight multiplicities and show that if the action of nonzero real vectors of  each affine component of the zero part is  neither completely injective nor completely locally nilpotent, then these modules  are parabolically induced.
\end{abstract}

\section{Introduction}
Representation theory  of Lie superalgebras is one of the most important topics which  mathematicians  and physicists are interested in.

The study of representations of a Lie superalgebra  $\LL=\LL_0\op\LL_1$ having root space decomposition with respect to a splitting Cartan subalgebra $\hh\sub \LL_0$, strongly depends on the existence of  parabolic subsets of the corresponding  root system $R$; i.e., those  subsets $P$  satisfying
 $R=P\cup -P$ and $(P+P)\cap R\sub P.$ {
For a proper parabolic subset} $P$ of $R$ and
 $$\hbox{\small $\LL^\circ:=\op_{\a\in P\cap -P}\LL^\a,\;\LL^+:=\op_{\a\in P\setminus-P}\LL^\a,\;\LL^-:=\op_{\a\in - P\setminus P}\LL^\a\andd \frak{p}:=\LL^\circ\op\LL^+,$}$$
  each $\LL^\circ$-module $N$ is a module over $\frak{p}$ with trivial action of $\LL^+.$ Then  $$\tilde N:=U(\LL)\ot_{U(\frak{p})}N$$ is an $\LL$-module; here  $U(\LL)$ and $U(\frak{p})$ denote respectively  the universal enveloping algebras of $\LL$ and $\frak{p}.$ If the $\LL$-module  $\tilde N$ contains a unique maximal submodule $Z$ intersecting $N$ trivially, the quotient module $${\rm Ind}_{\LL}(N):=\tilde N/Z$$ is called  a parabolically induced module.
{A cuspidal module is defined to be an irreducible $\LL$-module  which is not parabolically induced  from an irreducible module over $\LL^\circ$ corresponding to a parabolic subset of $R$}.
The first  key point in the study of irreducible representations of $\LL$ is { to recognize whether}  a module is cuspidal or not. In this case, the classification problem is reduced to the classification of cuspidal modules.
In this regard, the first concern after   classification of irreducible finite dimensional $\LL$-modules,  is the classification of finite weight $\LL$-modules, i.e., {those}  $\mathcal{L}$-modules $M$ having a weight space decomposition  $$M=\op_{\lam\in \mathcal{H}^\ast}M^\lam$$ in which $\mathcal{H}^\ast$ is the dual space of $\mathcal{H}$  and
$$M^\lam:=\{v\in M\mid hv=\lam(h)v\;\;(h\in\mathcal{H})\}\quad\quad(\lam\in \mathcal{H}^\ast)$$ is finite dimension.

{Irreducible finite weight  modules over  a {finite-dimensional}  reductive Lie algebra $\LL$ was studied by S.L. Fernando \cite{F} in 1990. He showed that the classification of such modules is reduced to the classification of irreducible finite weight modules on which all nonzero root vectors act nilpotently or all of them act injectively. He also proved that when all nonzero root vectors of a {finite-dimensional} simple Lie algebra $\LL$ act injectively on an irreducible finite weight module (i.e., the module is cuspidal), then $\LL$ is either of type $A$ or $C.$}

In   affine Lie algebra case, the existence of imaginary roots  makes the study  more complicated. Irreducible finite weight modules over affine Lie algebras have been {studied} in  \cite{C}--\cite{CP2}, \cite{E2}, \cite{F1}--\cite{F5}, \cite{FS} and \cite{DG}.

The super version of affine Lie algebras was introduced and classified by J.W.~Van de Leur in 1986.  Due to his classification, affine Lie superalgebras with nonzero odd part are  known as untwisted types  $X^{(1)},$ where $X$ is the type of a {finite-dimensional} basic classical simple Lie superalgebra with nonzero odd part or twisted types $A(2k-1,2\ell-1)^{(2)}$  ({\tiny$(k,\ell)\neq (1,1)$}), $A(2k,2\ell)^{(4)},$  $A(2k,2\ell-1)^{(2)},$ and $D(k,\ell)^{(2)}.$

{In 2003, S. Eswara Rao \cite{E1} proved that if the zero part of a basic classical simple Lie superalgebra $\fg$ has at least two simple components, then all integrable modules  over the untwisted affine Lie superalgebra based on $\fg$ with nonzero action of the canonical central element, are trivial.} Then in   2006,  S. Eswara Rao and V. Futorny \cite{F6}, \cite{EF} classified irreducible finite weight modules, over untwisted affine Lie superalgebras, on which the canonical central element acts as a nonzero multiple of the identity map. They reduced the classification problem to the classification of cuspidal modules of {finite-dimensional}  cuspidal {Lie superalgebras} which is discussed in a work of Dimitrov, Mathieu and Penkov \cite{DMP}.
  Recently, L. Calixto and  V. Futorny have studied highest weight modules over untwisted  affine Lie superalgebras \cite{CF}.
In \cite{you8}, we initiated the  study of finite weight modules  over  twisted affine Lie superalgebras.

An affine Lie superalgebra $\LL$ is equipped with a nondegenerate invariant supersymmetric bilinear form and has  a weight space decomposition with respect to a {finite-dimensional} abelian subalgebra $\hh\sub \LL_0$ on which the form is nondegenerate. This allows to transfer the form on $\LL$ to a nondegenerate bilinear form $\fm$ on the dual space $\hh^*$ of $\hh$  and  divide nonzero roots of the root system $R$ of $\LL$ (with respect to $\hh$) into three {parts:}
$R_{re}^\times$ (nonzero real roots), consisting of those roots $\a$ with $(\a,\a)\neq 0,$ $R_{im}^\times$ (nonzero imaginary roots), consisting of nonzero roots $\a$ with $(\a,R)=\{0\}$ and $R_{ns}^\times$ (nonzero nonsingular roots) consisting of all nonzero roots  which {are  neither real nor imaginary.}  Moreover, the set of imaginary roots of $\LL$ generates a free abelian group $\bbbz\d$ of rank 1.

The structure of an irreducible  finite weight module $M$  over the affine Lie superalgebra $\LL$ strongly depends on the nature of the action of root vectors corresponding to nonzero real roots. More precisely, each nonzero root vector corresponding to a nonzero real root $\a$, acts on $M$ either injectively or locally nilpotently. We denote by $R^{in}$ (resp. $R^{ln}$) the set of all nonzero real root $\a$ whose corresponding nonzero root vectors  act on $M$  injectively (resp.  locally nilpotently).  In \cite{you8}, we showed that  for each nonzero real root $\a,$ one of the following occurs:
\begin{itemize}
\item  $\a$ is full-locally nilpotent, i.e., $R\cap (\a+\bbbz\d)\sub R^{ln},$
\item $\a$ is full-injective, i.e.,  $R\cap (\a+\bbbz\d)\sub R^{in},$
\item  $\a$ is hybrid, i.e., it is neither full-locally nilpotent nor full-injective.
\end{itemize}
In particular, we can divide our study into two cases when all real roots are hybrid or not; we call an irreducible finite weight module over a twisted affine Lie superalgebra hybrid if all nonzero real roots are hybrid and otherwise call it tight.

{An irreducible finite weight module $V$ over} an affine Lie superalgebra $\LL$ is parabolically induced if
$$V^{\LL^+}:=\{v\in V\mid \LL^+ v=\{0\}\}$$ is nonzero for  $\LL^+=\op_{\a\in P\setminus-P}\LL^\a$ where $P$ is   certain parabolic subset; see Proposition~\ref{ind} for the details.

When the mentioned parabolic subset $P$  is determined, the most difficult part to show $V^{\LL^+}\neq\{0\}$ is dealing with nonsingular roots; more precisely, the parabolic subset is usually selected such that $(P\setminus-P)\cap R_{re}\sub R^{ln}$  and $(-P\setminus P)\cap R_{re}\sub R^{in}.$ We can  find nonzero weight vectors $v$  with $\LL^\a v=\{0\}$ for $\a\in (R_{re}\cup R_{im})\cap (P\setminus-P)$ and then we need  to prove that among these vectors, there are some $v$ annihilated by $\LL^\a$ for all $\a\in R_{ns}\cap (P\setminus-P).$ If $\d\in P\setminus-P,$ then
 for each nonsingular root $\a,$  we have $\a+m\d\in P\setminus-P$ for large enough $m.$ This is very helpful to handel the situation but this does not happen for tight modules which in turn makes the situation more complicated.

Another difficulty coming up in the study of tight irreducible  finite weight modules is that up to an $\hh$-module whose weights are  nonzero imaginary roots, the even part of a twisted affine Lie superalgebra $\LL\neq A(0,2\ell)^{(4)},A(0,2\ell-1)^{(2)},D(0,\ell)^{(2)}$  is a sum of two affine Lie algebras; say $\LL_0(1)$ and $\LL_0(2)$ with corresponding root systems $R_0(1)$ and $R_0(2)$ respectively,
and the study  depends on the fact that  if non-hybrid roots occur for roots of both $\LL_0(1)$ and $\LL_0(2)$ or just for the roots of one of them. If for some $i=1,2,$ all real  roots of $R_0(i)$ are hybrid, the structure of the modules again depends on whether  $R_0(j)\cap R^{ln}$ is a nonempty proper subset or not for $j$ with  $\{i,j\}=\{1,2\}.$

In \cite{you8}, we have studied hybrid irreducible finite weight modules over twisted affine Lie superalgebras $\LL$ with $\LL_1\neq \{0\}$ and here we focus on tight finite weight modules over $\LL.$
 We show that tight irreducible  finite weight $\LL$-modules with bounded weight multiplicities and proper nonempty subsets  $R_0(j)\cap R^{ln}\sub R_0(j)$ ($j=1,2$) are parabolically induced.
\section{Preliminary}
In this work, we deal with twisted affine Lie superalgebras with nonzero odd part, so in what follows we briefly explain their structures; any information stated without proof can be found in \cite{van-thes}; one can also see \cite[Appendix]{you8}. Suppose that  $\fg$ is a complex {finite-dimensional}  basic classical simple Lie superalgebra of type $X=A(k,\ell) (\hbox{\tiny$(k,\ell)\neq (1,1)$}), D(k,\ell)$ with a Cartan subalgebra $\fh\sub\fg_0;$ here $k$ is a nonnegative integer and $\ell$ is a positive integer. Suppose that $\kappa\fm$ is  a nondegenerate supersymmetric invariant  even  bilinear form on $\fg.$ In \cite{van-thes}, the author introduces
a certain automorphism $\sg:\fg\longrightarrow \fg$  such that
\begin{itemize}
\item{\small \rm $\sg$ is of order $n=4$ if $X=A(2k,2\ell),$}
\item{\small \rm $\sg$ is of order $n=2$ if $X=A(2k-1,2\ell-1),A(2k,2\ell-1),D(k,\ell).$}\\
\end{itemize}
Since $\sg$ preserves $\fg_0$ as well as $\fg_1,$ we have
\begin{align*}
\fg_i=\bigoplus_{j=0}^{n-1}{}^{[j]}\fg_i\quad\hbox{where}\quad {}^{[j]}\fg_i=\{x\in\fg_i\mid \sg(x)=\zeta^jx\}\quad\quad(i\in\bbbz_2,\; 0\leq j\leq n-1)
\end{align*}
in which $\zeta$ is the $n$-th primitive root of unity.  Then,
we have
\[\fg=\bigoplus_{j=0}^{n-1}{}^{[j]}\fg\quad\quad\hbox{ with } {}^{[j]}\fg={}^{[j]}\fg_0\op{}^{[j]}\fg_1\quad(0\leq j\leq n-1).\]
 Set $$\LL:=\bigoplus_{j=0}^{n-1}({}^{[j]}\fg\ot t^j\bbbc[t^{\pm n}])\op\bbbc c\op\bbbc d\andd \hh:=(({}^{[0]}\fg\cap\fh)\ot 1)\op\bbbc c\op\bbbc d, $$ then  $\LL,$ which is denoted by $X^{(n)},$ together with
$$[x\ot t^p+rc+sd,y\ot t^q+r'c+s'd]:=[x,y]\ot t^{p+q}+p\kappa(x,y)\d_{p+q,0}c+sqy\ot t^q-s'px\ot t^p$$ is called  a twisted affine Lie superalgebra of type $X^{(n)}.$ The Lie superalgebra $\LL$ has a  weight space decomposition with respect to  $\hh.$ We denote  the corresponding root system by $R.$

The form $\kappa\fm$ induces the following nondegenerate  supersymmetric invariant  bilinear form  $\fm$  on $\LL:$
$$(x\ot t^p+rc+sd,y\ot t^q+r'c+s'd):=\kappa(x,y)\d_{p+q,0}+rs'+sr'.$$
 As the form is nondegenerate on  $\hh,$ one can transfer the  form  on  $\hh$ to a form on $\hh^*$ denoted again by $\fm.$
   Setting
\begin{equation}\label{im-ns}
\begin{array}{lll}
R_{re}^\times := \{\a\in R\mid (\a,\a)\neq 0\},& R_{re}:=\{0\}\cup R_{re}^\times & \hbox{(real roots)},\\
R_{im}:=\{\a\in R\mid (\a,\b)=0\;\;\forall \b\in R\},& R_{im}^\times:=R_{im}\setminus\{0\}& \hbox{(imaginary  roots)},\\
R_{ns}:=\{0\}\cup (R\setminus (R_{re}\cup R_{im})),& R_{ns}^\times:=R_{ns}\setminus\{0\}& \hbox{(nonsingular  roots)},
\end{array}
\end{equation}
we have
$R=R_{im}\cup R_{re}\cup R_{ns}.$ It is known that  $R_{im}$ generates a free abelian group of rank 1; say $\bbbz\d.$ Also,
\begin{equation}
\label{dim}
\dim(\LL^\a)=1\quad\quad(\a\in R\setminus R_{im})
\end{equation}
and
\begin{equation}
\label{sl2}
\parbox{2.8in}{if $\a\in R_{re}^\times\cap R_0,$ then there are $e\in\LL^\a$ and $f\in\LL^{-\a}$ such that $(e,f,[e,f])$ is an $\frak{sl}_2$-triple.}
\end{equation}

Moreover, the root system $R$ of $X^{(n)}$  has an expression  as   in the following table:

\newpage
  \begin{table}[h]\caption{} \label{table1}
 {\footnotesize \begin{tabular}{|c|l|}
\hline
$X^{(n)}$ &\hspace{3.25cm}$R$ \\
\hline
$A(2k,2\ell-1)^{(2)}$&$\begin{array}{rcl}
\bbbz\d
&\cup& \bbbz\d\pm\{\ep_i,\d_j,\ep_i\pm\ep_r,\d_j\pm\d_s,\ep_i\pm\d_j\mid i\neq r,j\neq s\}\\
&\cup& (2\bbbz+1)\d\pm\{2\ep_i\mid 1\leq i\leq k\}\\
&\cup& 2\bbbz\d\pm\{2\d_j\mid 1\leq j\leq \ell\}.
\end{array}$\\
\hline
$A(2k-1,2\ell-1)^{(2)},\;(k,\ell)\neq (1,1)$& $\begin{array}{rcl}
\bbbz\d&\cup& \bbbz\d\pm\{\ep_i\pm\ep_r,\d_j\pm\d_s,\d_j\pm\ep_i\mid i\neq r,j\neq s\}\\
&\cup& (2\bbbz+1)\d\pm\{2\ep_i\mid 1\leq i\leq k\}\\
&\cup& 2\bbbz\d\pm\{2\d_j\mid 1\leq j\leq \ell\}
\end{array}$\\
\hline
$A(2k,2\ell)^{(4)}$& $\begin{array}{rcl}
\bbbz\d&\cup&  \bbbz\d\pm\{\ep_i,\d_j\mid 1\leq i\leq k,\;1\leq j\leq \ell\}\\
&\cup& 2\bbbz\d\pm\{\ep_i\pm\ep_r,\d_j\pm\d_s,\d_j\pm\ep_i\mid i\neq r,j\neq s\}\\
&\cup&(4\bbbz+2)\d\pm\{2\ep_i\mid 1\leq i\leq k\}\\
&\cup& 4\bbbz\d\pm\{2\d_j\mid 1\leq j\leq \ell\}
\end{array}$\\
\hline
$D(k+1,\ell)^{(2)}$& $\begin{array}{rcl}
\bbbz\d&\cup&  \bbbz\d\pm\{\ep_i,\d_j\mid 1\leq i\leq k,\;1\leq j\leq \ell\}\\
&\cup& 2\bbbz\d\pm\{2\d_j,\ep_i\pm\ep_r,\d_j\pm\d_s,\d_j\pm\ep_i\mid i\neq r,j\neq s\}
\end{array}$\\
\hline
 \end{tabular}}
 \end{table}
with
\[R_{ns}^\times=R\cap (\bbbz\d\pm\{\ep_i\pm\d_j\mid 1\leq i\leq k,\;1\leq j\leq \ell\}).\] One can see that
\begin{equation}\label{re-re}
(R^\times_{ns}+R^\times_{ns})\cap R\sub R_{re}\cup R_{im}.
\end{equation}

We also have from Table~\ref{table1} that  $R\sub \dot R+\bbbz\d$ where $\dot R$ is as in the following table:

{\tiny
\begin{table}[h]
{\footnotesize\begin{tabular}{|c|c|}
\hline
$X^{(n)}$ & $\dot R$ \\
\hline
$A(2k,2\ell-1)^{(2)}$& $\pm\{\ep_i,\d_j,\ep_i\pm\ep_r,\d_j\pm\d_s,\ep_i\pm\d_j\mid 1\leq i,r\leq k,\; 1\leq j,s\leq \ell\}$\\
\hline
$A(2k-1,2\ell-1)^{(2)}$& $\pm\{\ep_i\pm\ep_r,\d_j\pm\d_s,\ep_i\pm\d_j\mid 1\leq i,r\leq k,\; 1\leq j,s\leq \ell\}$\\
$(k,\ell)\neq (1,1)$&\\
\hline
$A(2k,2\ell)^{(4)}$& $\pm\{\ep_i,\d_j,\ep_i\pm\ep_r,\d_j\pm\d_s,\ep_i\pm\d_j\mid 1\leq i,r\leq k,\; 1\leq j,s\leq \ell\}$\\
\hline
$D(k+1,\ell)^{(2)}$&$\pm\{\ep_i,\d_j,\ep_i\pm\ep_r,\d_j\pm\d_s,\ep_i\pm\d_j\mid 1\leq i\neq r\leq k,\; 1\leq j,s\leq \ell\}$\\
\hline
 \end{tabular}
}
 \end{table}
 }

A nonzero root $\dot\a\in\dot R$ is called a nonzero real (resp. nonsingular, imaginary) root if $(\dot\a+\bbbz\d)\cap R\sub R_{re}$ (resp. $R_{ns}, R_{im}$). Setting
\begin{equation}\label{s}S_{\dot\a}:=\{\sg\in\bbbz\d\mid \dot\a+\sg\in R\}\quad\quad(\dot\a\in \dot R),\end{equation} we get

\begin{table}[h]\caption{} \label{table2}
\footnotesize{\begin{tabular}{|l|c|c|c|c|}
  \hline
  & $A(2k,2\ell-1)^{(2)}$& $A(2k-1,2\ell-1)^{(2)}$&$A(2k,2\ell)^{(4)}$&$D(k+1,\ell)^{(2)}$\\
  \hline
  $S_{\pm\ep_i}$&$\bbbz\d$&$\emptyset$&$\bbbz\d$&$\bbbz\d$\\
  \hline
  $S_{\pm\ep_i\pm\ep_j}$&$\bbbz\d$&$\bbbz\d$&$2\bbbz\d$&$2\bbbz\d$\\
  \hline
  $S_{\pm2\ep_i}$&$(2\bbbz+1)\d$&$(2\bbbz+1)\d$&$(4\bbbz+2)\d$&$\emptyset$\\
  \hline
  $S_{\pm\d_j}$&$\bbbz\d$&$\emptyset$&$\bbbz\d$&$\bbbz\d$\\
  \hline
  $S_{\pm\d_j\pm\d_q}$&$\bbbz\d$&$\bbbz\d$&$2\bbbz\d$&$2\bbbz\d$\\
  \hline
  $S_{\pm2\d_j}$&$2\bbbz\d$&$2\bbbz\d$&$4\bbbz\d$&$2\bbbz\d$\\
  \hline
  $S_{\pm\ep_i\pm\d_j}$&$\bbbz\d$&$\bbbz\d$&$2\bbbz\d$&$2\bbbz\d$\\
  \hline
\end{tabular}}
\end{table}
One can see that for each $\dot\a\in \dot R\setminus\{0\},$ there is $r_{\dot\a}\in\{1,2,4\}$
 and $0\leq k_{\dot\a}<r_{\dot\a}$ such that
 \begin{equation}\label{salpha}S_{\dot\a}=(r_{\dot\a}\bbbz+k_{\dot\a})\d.\end{equation}

We conclude this section with some information regarding  $\LL_0$ and its root system $R_0.$
Up to an $\hh$-submodule of $\LL$ for which all weights are nonzero imaginary roots,   $\LL_0$ is a sum of two affine Lie algebras $\LL_0(1)$ and $\LL_0(2)$ with  $\LL_0(1)\cap \LL_0(2)=\bbbc c+\bbbc d$ if $k\neq 0$ and  an
affine Lie algebra if $k=0$; see \cite[Appendix]{you8} for the details. We  have $R_0=R_0(1)\cup R_0(2)$  in which $R_0(2)=\emptyset$ if $k=0$ and  $R_0(1)$  (the root system of $\LL_0(1))$ and $R_0(2)$  (the root system of $\LL_0(2)$ if  $k\neq 0)$  are given as in the following table:

 \begin{table}[h]\caption{}\label{table3}
{\footnotesize \begin{tabular}{|c|l|l|}
\hline
$X^{(n)}$& \multicolumn{1}{>{\centering\arraybackslash}m{45mm}|}{$R_0(1)$}  &\multicolumn{1}{>{\centering\arraybackslash}m{45mm}|}{$R_0(2)\hbox{ if $k\neq 0$}$}\\
\hline
$A(2k,2\ell-1)^{(2)}$
& $\begin{array}{l}
(2\d_{\ell,1}+(1-\d_{\ell,1}))\bbbz\d\\
\cup \bbbz\d\pm\{\d_j\pm\d_s\mid 1\leq j\neq s\leq \ell\}\\
\cup 2\bbbz\d\pm\{2\d_j\mid 1\leq j\leq \ell\}
\end{array}$
&$\begin{array}{l}
\bbbz\d\\
\cup \bbbz\d\pm\{\ep_i,\ep_i\pm\ep_r\mid 1\leq i\neq r\leq k\}\\
\cup (2\bbbz+1)\d\pm\{2\ep_i\mid 1\leq i\leq k\}
\end{array}$
\\
\hline
$\begin{array}{c}
A(2k-1,2\ell-1)^{(2)}\\(k,\ell)\neq (1,1)
\end{array}
$
& $\begin{array}{l}
(2\d_{\ell,1}+(1-\d_{\ell,1}))\bbbz\d \\
\cup \bbbz\d\pm\{\d_j\pm\d_s\mid 1\leq j\neq s\leq \ell\}\\
\cup 2\bbbz\d\pm\{2\d_j\mid 1\leq j\leq \ell\}
\end{array}$
& $\begin{array}{l}
(2\d_{k,1}+(1-\d_{k,1}))\bbbz\d\\
\cup \bbbz\d\pm\{\ep_i\pm\ep_r\mid 1\leq i\neq r\leq k\}\\
\cup (2\bbbz+1)\d\pm\{2\ep_i\mid 1\leq i\leq k\}
\end{array}$
\\
\hline
$A(2k,2\ell)^{(4)}$
& $\begin{array}{l}
2\bbbz\d\\
\cup (2\bbbz+1)\d\pm\{\d_j\mid 1\leq j\leq \ell\}\\
\cup 2\bbbz\d\pm\{\d_j\pm\d_s\mid 1\leq j\neq s\leq \ell\}\\
\cup 4\bbbz\d\pm\{2\d_j\mid 1\leq j\leq \ell\}
\end{array}$
& $\begin{array}{l}
2\bbbz\d\\
\cup  2\bbbz\d\pm\{\ep_i\mid 1\leq i\leq k\}\\
\cup 2\bbbz\d\pm\{\ep_i\pm\ep_r\mid 1\leq i\neq r\leq k\}\\
\cup(4\bbbz+2)\d\pm\{2\ep_i\mid 1\leq i\leq k\}
\end{array}$
\\
\hline
$
D(k+1,\ell)^{(2)}
$
& $\begin{array}{l}
2\bbbz\d\\
\cup 2\bbbz\d\pm\{\d_j\pm\d_s\mid 1\leq j, s\leq \ell\}
\end{array}$
& $\begin{array}{l}
\bbbz\d\\
\cup \bbbz\d\pm\{\ep_i\mid {1\leq i\leq k}\}\\
\cup 2\bbbz\d\pm\{\ep_i\pm\ep_r\mid 1\leq i\neq r\leq k\}
\end{array}$
\\
\hline
 \end{tabular}}
 \end{table}
\noindent We  see that
\begin{equation}
\label{span}
\sspan_\bbbr R_0=\sspan_\bbbr R=\sspan_\bbbr\{\d,\ep_i,\d_j\mid 1\leq i\leq k,1\leq j\leq \ell\}.
\end{equation}
Also  setting
\[\dot R_0(i):=\{\dot \a\in \dot R\mid \exists \sg\in \bbbz\d\ni \dot\a+\sg\in R_0(i)\}\quad\quad(i=1,2),\] we have

\begin{table}[h]\caption{}\label{table4}
{\footnotesize \begin{tabular}{|c|c|c|}
\hline
$X^{(n)}$ & $\dot R_0(1)$ &$\dot R_0(2)$\\
\hline
$A(2k,2\ell-1)^{(2)}$
& $
 \{\pm\d_j\pm\d_s\mid 1\leq j, s\leq \ell\}
$&$
\{\pm\ep_i,\pm\ep_i\pm\ep_r\mid 1\leq i,r\leq k\}
$
\\
\hline
$\begin{array}{c}
A(2k-1,2\ell-1)^{(2)}\\(k,\ell)\neq (1,1)
\end{array}
$
& $
\pm\{\d_j\pm\d_s\mid 1\leq j, s\leq \ell\}$& $ \{\pm\ep_i\pm\ep_r\mid 1\leq i,r\leq k\}$
\\
\hline
$A(2k,2\ell)^{(4)}$
& $\{\pm\d_j,\pm\d_j\pm\d_s\mid 1\leq j, s\leq \ell\}$& $
 \{\pm\ep_i,\pm\ep_i\pm\ep_r\mid 1\leq i, r\leq k\}
$
\\
\hline
$
D(k+1,\ell)^{(2)}
$
& $\{\pm\d_j\pm\d_s\mid 1\leq j, s\leq \ell\}
$& $\{\pm\ep_i,\pm\ep_i\pm\ep_r\mid 1\leq i\neq r\leq k\}
$
\\
\hline
 \end{tabular}}
 \end{table}

For $i=1,2,$ $\dot R_0(i)$ is an irreducible finite root system provided that it is nonempty. Setting
\begin{equation}
\label{s-i}
S_{\dot\a}(i):=\{\sg\in\bbbz\d\mid \dot\a+\sg\in R_0(i)\}\quad\quad(\dot\a\in \dot R_0(i);\;\;i=1,2),
\end{equation} we have

\newpage

\begin{table}[h]\caption{}\label{table5}
\footnotesize{\begin{tabular}{|l|c|c|c|c|}
  \hline
  & $A(2k,2\ell-1)^{(2)}$& $A(2k-1,2\ell-1)^{(2)}$&$A(2k,2\ell)^{(4)}$&$D(k+1,\ell)^{(2)}$\\
  \hline
  $S_{\pm\d_p}(1)$&$\emptyset$&$\emptyset$&$(2\bbbz+1)\d$&$\emptyset$\\
  \hline
  $S_{\pm\d_p\pm\d_q}(1)$&$\bbbz\d$&$\bbbz\d$&$2\bbbz\d$&$2\bbbz\d$\\
  \hline
  $S_{\pm2\d_p}(1)$&$2\bbbz\d$&$2\bbbz\d$&$4\bbbz\d$&$2\bbbz\d$\\
  \hline
  $S_{\pm\ep_i}(2)$&$\bbbz\d$&$\emptyset$&$2\bbbz\d$&$\bbbz\d$\\
  \hline
  $S_{\pm\ep_i\pm\ep_j}(2)$&$\bbbz\d$&$\bbbz\d$&$2\bbbz\d$&$2\bbbz\d$\\
  \hline
  $S_{\pm2\ep_i}(2)$&$(2\bbbz+1)\d$&$(2\bbbz+1)\d$&$(4\bbbz+2)\d$&$\emptyset$\\
  \hline
\end{tabular}}
\end{table}
Denoting the root length of a root $\dot\a$ of an irreducible finite root system by $\ell(\dot\a),$ it is readily seen from this table that
\begin{equation}
\label{sum}
\hbox{\parbox{4in}{\begin{center}if $\dot\a,\dot\b,\dot\a+\dot\b\in\dot R_0(i)\setminus\{0\}~ (i=1,2)$ with $\ell(\dot\a)=\ell(\dot\b)\leq \ell(\dot\a+\dot\b),$ then $S_{\dot\a+\dot\b}(i)\sub S_{\dot\a}(i)+S_{\dot\b}(i).$\end{center} } }
\end{equation}
\begin{rem}\label{rem}
{\rm For our further use, we point out the following:
\begin{itemize}
\item[(i)] For types $A(2k,2\ell-1)^{(2)},$ $A(2k-1,2\ell-1)^{(2)}$ and $A(2k,2\ell)^{(4)},$ recalling (\ref{salpha}), we set
 \[r:=r_{\pm\ep_i\pm\d_p}\quad\quad(1\leq i\leq k,\; 1\leq p\leq \ell).\] Then from Tables~\ref{table2},~\ref{table5}, we have
 \[S_{\pm\ep_i\pm\d_p}=r\bbbz\d,\;\;
 S_{\pm2\ep_i}(2)=(2r\bbbz+r)\d\andd S_{\pm2\d_p}(1)=2r\bbbz\d.\] This together with the fact that  each nonzero nonsingular root of $\dot R$ is of the form $\pm\ep_i\pm\d_p,$ implies that if $\dot\eta\in \dot R_{ns}\setminus\{0\},$ then   there are $\dot\a,\dot \b$  with
 \[\pm 2\dot\a+2r_{\dot\eta}\bbbz\d\sub R_0\cap R_{re},\;\;\; \pm2\dot\b\pm r_{\dot\eta}\d+2r_{\dot\eta}\bbbz\d\sub R_0\cap R_{re},\]
 \[\dot\eta=\dot\a+\dot\b\andd 2\dot\a\pm 2\dot\b\not \in \dot R.\]
 \item[(ii)]
 If $R=D(k+1,\ell)^{(2)},$  then for $1\leq i\leq k$ and $1\leq p\leq \ell,$ we have $r_{\pm\ep_i\pm\d_p}=2,$ more precisely, we have
 \[S_{\pm\ep_i\pm\d_p}=2\bbbz\d,\;\;S_{\pm\ep_i}(2)=\bbbz\d\andd S_{\pm2\d_p}(1)=2\bbbz\d.\]
So if $\dot\eta\in \dot R_{ns}^\times,$ there are $1\leq i\leq k$ and $1\leq p\leq \ell$ such that for $\dot\a\in\{\pm\ep_i\}$ and $\dot\b\in\{\pm\d_p\},$  we have $\dot\eta=\dot\a+\dot\b,$     \[\pm\dot\a+\underbrace{r_{\dot\eta}\bbbz\d}_{\sub\bbbz\d}\sub R_0\cap R_{re},\;\;\; \pm2\dot\b\underbrace{\pm r_{\dot\eta}\d+2r_{\dot\eta}\bbbz\d}_{\sub2\bbbz\d}\sub R_0\cap R_{re}\andd \dot\a\pm 2\dot\b\not\in \dot R.\]
 \end{itemize}}
 \end{rem}

\section{Finite weight modules having shadow}\label{generic}
Throughout this section, we assume  $\LL=\LL_0\op\LL_1$ is a twisted  affine Lie superalgebra  with   Cartan subalgebra $\hh\sub \LL_0$ such that $\LL_1\neq \{0\}$ {and the root system   from Table~\ref{table1}.} We keep the notations
introduced in the previous section.

A module $M$ over $\LL$ is called a {\it weight module} if it is equipped with  a weight space decomposition with respect to $\hh,$ that is,
$M=\op_{\lam\in\hh^*}M^\lam$ where for $\lam$ in the dual space $\hh^*$ of $\hh,$
\[M^\lam:=\{v\in M\mid hv=\lam(h)v\;\;\;(\forall h\in\hh)\}.\]
 Each element of the {\it support}
 \[\supp(M):=\{\lam\in\hh^*\mid M^\lam\neq\{0\}\}\] of $M$  is called a  {\it weight} of $M.$ Each $M^\lam$ is called a {\it weight space} and each vector of a weight space is called a {\it weight vector}.
  The module $M$ is called a {\it finite weight module} if each weight space is {finite-dimensional}.

  {Assume} $M$ is a weight $\LL$-module. Denote by $R^{in}$ (resp. $R^{ln}$) the set of all nonzero $\a\in R_{re}$ for which  $0\neq x\in \LL^\a$ acts injectively (resp. locally nilpotently) on $M.$
It is easily checked that if $M$ is irreducible, then $R_{re}^\times=R^{in}\cup R^{ln}.$

\begin{deft}\label{gen}
{\rm
\begin{itemize}
\item[(i)]
Suppose that $S\sub R.$ We say a decomposition $S=S^+\cup S^\circ\cup S^-$ is a {\it triangular} decomposition for $S$ if there is a linear functional $\boldsymbol\zeta:\hbox{span}_\bbbr S\longrightarrow \bbbr$ such that {\small $$S^{+}=\{\a\in S\mid \boldsymbol\zeta(\a)>0\},\;S^{-}=\{\a\in S\mid \boldsymbol\zeta(\a)<0\}\andd S^\circ=\{\a\in S\mid \boldsymbol\zeta(\a)=0\}.$$}
The decomposition is called {\it trivial} if $S=S^\circ.$
\item[(ii)]  Assume $M$ is a weight  $\LL$-module.
Set

\noindent {\small$\frak{B}_M:=$}{\small\hbox{$\{\a\in \hbox{span}_\bbbz R\mid  \{k\in\bbbz^{>0}\mid \lam+k\a\in \supp(M)\}$ is finite for all $\lam\in\supp(M)$\}}}

and

\noindent \noindent{\small$\frak{C}_M:=$}{\small$\{\a\in \hbox{span}_\bbbz R\mid   \a+\supp(M)\sub\supp(M)\}.$}

\smallskip

 We say  $M$ has {\it shadow} if
 \begin{itemize}
 \item[\rm\bf (s1)]
$R_{re}\setminus\{0\}=R^{in}\cup R^{ln},$
\item[\rm\bf (s2)] $R^{ln}=\frak{B}_M\cap R_{re}^\times$ and $R^{in}=\frak{C}_M\cap R_{re}^\times.$
 \end{itemize}
 \end{itemize}}
 \end{deft}

\begin{Pro}[{\cite[Pro.~3.3]{you8}}]\label{ind}
 Suppose that $R=R^+\cup R^\circ\cup R^-$ is a nontrivial  triangular decomposition for $R$ and $R^\circ=R^{\circ,+}\cup R^{\circ,\circ}\cup R^{\circ,-}$ is a triangular decomposition for $R^\circ.$ Set
$$\LL^\circ=\op_{\a\in R^{\circ,\circ}}\LL^\a,\;\; \LL^\pm=\op_{\a\in R^\pm\cup R^{\circ,\pm}}\LL^\a\andd \frak{p}=\LL^\circ\op\LL^+.$$
\begin{itemize}
\item[(i)] If $N$ is a nonzero  weight module over $\LL^\circ$ whose support lies in a single coset of $\sspan_\bbbz R^{\circ,\circ},$ then
 $$\tilde N:=U(\LL)\ot_{U(\frak{p})}N$$ has a unique maximal submodule $Z$ intersecting $N$ trivially. Moreover, the induced module  $${\rm Ind}_\LL(N):=\tilde{N}/Z$$ is an irreducible   $\LL$-module if and only if $N$ is an irreducible $\LL^\circ$-module.
\item[(ii)] If $V$ is an irreducible finite weight $\LL$-module with $$V^{\LL^+}:=\{v\in V\mid \LL^+ v=\{0\}\}\neq \{0\},$$ then $V^{\LL^+}$ is an irreducible finite weight $\LL^\circ$-module and  $V\simeq {\rm Ind}_\LL(V^{\LL^+}).$
 \end{itemize}
\end{Pro}

\begin{Pro}\label{pdelta}
  Suppose that $R=R^+\cup R^\circ\cup R^-$ is  a triangular decomposition for $R$ and set
     $$R_{re}^{\pm}:=R^\pm\cap R_{re},\; R_{ns}^{\pm}:=R^\pm\cap R_{ns}\andd R_{im}^\pm:=R_{im}\cap R^\pm.$$
     Assume  $M$ is a weight module with $R_{re}^+ \sub \frak{B}_M$ and $R^-_{re}\sub \frak{C}_M.$
If  $\d\in R^\circ,$ then there is $\mu\in\supp(M)$ such that $(\mu+R_{re}^+)\cap\supp(M)=\emptyset.$
\end{Pro}
\pf For $\dot\a\in \dot R\setminus\{0\},$ suppose $r_{\dot\a}$ is as in (\ref{salpha}) and set
\begin{equation}\label{r}
r:=max\{r_{\dot\a}\mid \dot\a\in \dot R\setminus\{0\}\}.
\end{equation}
We know  from Table~\ref{table2}  that for each $\dot\a\in \dot R\setminus\{0\},$
$r_{\dot\a}|r$ and
$$\{m\in\bbbz\mid \dot\a+m\d\in R\}=\cup_{i=1}^{n_{\dot\a}}(r\bbbz+k_{\dot\a}^i)$$ for some $1\leq n_{\dot\a}\leq r$
 and $k_{\dot\a}^i\in\{0,\ldots,r-1\}$ $(1\leq i\leq n_{\dot\a}).$ Set $$\dot\Phi:=\{\dot\a+k_{\dot\a}^i\d\mid \dot\a\in \dot R_{re}\setminus\{0\}, 1\leq i\leq n_{\dot\a}\}\sub R.$$ Then, we have
 \begin{equation}\label{equal}
 R^\times=R\setminus R_{im}=\bigcup_{\dot\a\in\dot R\setminus\{0\} }\bigcup_{i=1}^{n_{\dot\a}}(\dot\a+k_{\dot\a}^i\d+r\bbbz\d)=\bigcup_{\dot\a\in\dot\Phi}(\dot\a+r\bbbz\d).
 \end{equation}

\noindent{\bf Claim 1.} Set $\dot\Phi^+:=\dot\Phi\cap R_{re}^+$ and   $$\Pi:=\{\dot\a\in\dot \Phi^+\mid \nexists \dot\b,\dot\gamma\in\dot \Phi^+\ni \dot\a=\dot\b+\dot\gamma\}.$$ Then for each element $\dot\b\in \dot \Phi^+,$ there are $t_{\dot\a}\in\bbbz^{\geq 0}$ ($\dot\a\in \Pi$) with $\dot\b=\sum_{\dot\a\in \Pi}t_{\dot\a}\dot\a:$ Suppose that $\boldsymbol\zeta$ is the function defining the triangular decomposition $R=R^+\cup R^\circ\cup R^-.$
\begin{equation}\label{less}
\parbox{4in}{ For $\dot\a,\dot\b,\dot\gamma\in \dot \Phi^+,$ if $\dot\b=\dot\a+\dot\gamma,$ then $\boldsymbol\zeta(\dot\a), \boldsymbol\zeta(\dot\gamma)\lneq\boldsymbol\zeta(\dot\b).$}
\end{equation}

Since $\dot\Phi$ is finite, we pick $\zeta_1<\cdots<\zeta_n$ such that $\{\boldsymbol\zeta(\dot\a)\mid \dot\a\in \dot \Phi^+ \}=\{\zeta_1,\ldots,\zeta_n\}.$
Suppose  $\dot \b\in \dot\Phi^+.$  If $\boldsymbol\zeta(\dot\b)=\zeta_1,$ (\ref{less}) implies  that $\dot\b\in \Pi$ and so we are done.  Next suppose $\dot\b\in \dot \Phi^+\setminus\Pi,$ then there are $\dot\a,\dot\gamma\in \dot \Phi^+$ with $\dot\b=\dot\a+\dot\gamma$ and so (\ref{less}) together with an induction process completes the proof of this step.

\noindent{\bf Claim  2.} There is $\lam\in \supp(M)$ such that $(\lam+\hbox{span}_{\bbbz^{\geq0}}\Pi)\cap \supp(M)=\{\lam\}:$ Use \cite[Pro.~3.6(ii)]{you8}.

\noindent{\bf Claim 3.} Consider $\lam$ as in Claim 2. If $\lam+\a\not\in \supp(M)$ for all $\a\in R_{re}^+,$ we get the result, otherwise, there is $\a\in R_{re}^+$ such that $\theta:=\lam+\a\in \supp(M).$ We claim that $\theta +\b\not \in \supp(M)$ for all $\b\in R_{re}^+.$ Suppose that $\b\in R_{re}^+,$ we  shall show
$\theta+\b\not\in \supp(M).$ Since $\a,\b\in R_{re}^+$ and $\d\in R^\circ,$ (\ref{equal}) implies that there are $\dot\a,\dot\b\in \dot \Phi^+$ and $m,n\in \bbbz$ such that
\[\a=\dot\a+mr\d\andd \b=\dot\b+nr\d.\]
   If to the contrary, $\theta+\b\in \supp(M),$ then $\lam+\dot\a+\dot\b+r(m+n)\d\in \supp(M).$ But $-(\dot\b+r(m+n)\d)\in R_{re}^-\sub \frak{C}_M,$ so $$\lam+\dot\a=(\underbrace{\lam+\dot\a+\dot\b+r(m+n)\d}_{\in \supp(M)})-(\dot\b+r(m+n)\d)\in (\lam+\dot \Phi^+)\cap \supp(M)$$ which gives a contradiction regarding Steps 1,2. This completes the proof.
\qed

\medskip

{\small \textbf{From now on till the end of this section, we assume $\boldsymbol{M}$ is a finite weight module over $\boldsymbol{\LL}$ having shadow.}}

 \begin{Thm}\label{BM}  Suppose $\a\in R_0\cap R^{ln},$ $-\a\in R^{in}$ and $\b_1,\ldots,\b_t\in R^{in}$ with $(\a,\b_i)=0$ and  $\a\pm\b_i\not \in R$ for $1\leq i\leq t.$
If there is  $\lam\in\supp(M)$ and a sequence
\begin{equation*}
0<s_1<s_2<\cdots
\end{equation*}
of positive integers such that $\{\lam+s_n\a+(s_n+1)(\b_1+\cdots+\b_t)\mid n\in\bbbz^{\geq0}\}\sub \supp(M),$ then dimensions of weight spaces are not bounded.
\end{Thm}
\pf
Suppose   $\a$ and $\b_1,\ldots,\b_t$ as well as $\lam\in \supp(M)$ and
\begin{equation}
\label{si}
0<s_1<s_2<\cdots
\end{equation} are as in the statement and set
\begin{equation}\label{b}
\b:=\b_1+\cdots+\b_t.
\end{equation}
Since $\a\in R_0,$ by (\ref{sl2}), there are $e\in\LL^\a$ and $f\in\LL^{-\a}$ such that $(e,f,h:=[e,f])$ is an $\frak{sl}_2$-triple. Also as  $\a\in R^{ln}$ and $-\a\in R^{in},$
\begin{equation}
\label{nil-inj}
\hbox{$e$ acts on $M$ locally nilpotently while $f$ acts injectively.}
\end{equation}
Set  $\fg:=\bbbc e\op \bbbc h\op \bbbc f$ and recall $\b$ from (\ref{b}), then
\begin{equation}\label{mn}M_n:=\bigoplus_{m\in\bbbz}M^{\lam+m\a+(s_n+1)\b}
\end{equation}  is a $\fg$-module  having a weight space decomposition with respect to $\bbbc h.$ As by our assumption,  $(\a,\b)=0,$ the set of weights of $K_n$ with respect to $\bbbc h$ is $$\{\lam(h)+2m\mid m\in\bbbz,\;\lam+m\a+(s_n+1)\b\in \supp(M)\}.$$ We mention that
\begin{equation}\label{finite}
\parbox{3.9in}{the weight space  corresponding to the weight $\lam(h)+2m$ ($m\in\bbbz$) is $M^{\lam+m\a+(s_n+1)\b}$
which is {finite-dimensional} by our assumption.}
\end{equation}

\medskip

\noindent{\bf Claim~1.} For each positive integer $n,$ there is a positive integer $m_n\geq s_n$ which is the largest positive integer with $\lam+m_n\a+(s_n+1)\b\in \supp(M).$

\smallskip

\noindent{\it Reason:}
We recall that  $M$ has shadow, so as   $\lam+s_n\a+(s_n+1)\b\in\supp(M)$ and $\a\in R^{ln}\sub\frak{B}_M,$   $$\{m\in\bbbz^{>0}\mid \lam+s_n\a+(s_n+1)\b+m\a\in \supp(M)\}$$ is finite. If for infinitely many positive integers $m,$ $\lam+m\a+(s_n+1)\b\in \supp(M),$ there are infinitely many positive integers $m$ with $\lam+s_n\a+(s_n+1)\b+m\a\in \supp(M)$ which is a contradiction.
Setting    $m_n$ to be the largest positive integer with $\lam+m_n\a+(s_n+1)\b\in \supp(M),$ since $\lam+s_n\a+(s_n+1)\b\in\supp(M),$ we have $s_n\leq m_n.$

\medskip

\noindent{\bf Claim~2.} Suppose $n$ is a positive integer, recall (\ref{mn}) and  assume $m_n$ is as in Claim~1. Set $$W(n):=\hbox{$\fg$-submodule of $M_n$ generated by $M^{\lam+m_n\a+(s_n+1)\b}.$}$$ Then
\[\hbox{$\lam(h)+2m_n$ is a weight for $W(n)$}\] and   $$W(n)=\op_{i\in I_n} W^{(\tau(i,n))}$$ where $I_n$ is an index set,  $\tau(i,n)\in\bbbc$ ($i\in I_n$) { is an eigenvalue  of the Casimir element} $c:=(h+1)^2+4fe$ of $\fg$ on $W(n)$ and $W^{(\tau(i,n))}$ is a $\fg$-submodule of {$W(n)$} whose weight vectors are generalized eigenvectors of $c$ corresponding to the eigenvalue $\tau(i,n).$
\smallskip

\noindent{\it Reason:} The first assertion is trivial and the last one follows from representation theory of $\frak{sl}_2$ (see \cite[\S 3.6]{Maz}) together with the fact  that $W(n)$ is a finitely generated  $\fg$-module having a weight space decomposition with respect to $\bbbc h$ with the set of weights $$\Lam_n\sub \lam(h)+2\bbbz$$ whose weight spaces are {finite-dimensional}; \cite[Pro.~2.1.1]{MP}.

\medskip

\noindent{\bf Claim~3.} Recall $I_n$ from Claim~2 and for each $n,$ {choose $k_n\in I_n$  such that} $\lam(h)+2m_n$ is a weight for $W^{(\tau(k_n,n))}.$ Then
\begin{equation}\label{new}
M^{\lam+m_{n}\a+(s_{n}+1)\b}\cap W^{(\tau(k_{n},n))}\neq \{0\}\quad\quad(n\in\bbbz^{>0})
\end{equation} and
$$
(\lam(h)+2m_n+1)^2= \tau(k_n,n).$$
\smallskip

\noindent{\it Reason:} The first assertion is obvious due to (\ref{finite}). For the last assertion, using  \cite[Pro. 3.55(v)]{Maz}, we get that   $W^{(\tau(k_n,n))}$ has finite length, so there are a positive integer $k$ and  submodules $W_i(n)$ ($1\leq i\leq k$) such that $$\{0\}:=W_0(n)\subsetneq W_1(n)\subsetneq W_2(n)\subsetneq \cdots\subsetneq W_k(n)=W^{(\tau(k_n,n))}$$ and $W_i(n)/W_{i-1}(n)$ is irreducible for all $1\leq i\leq k.$

If $\tau({k_n},n)\neq (\mu+1)^2$ for all $\mu \in \lam(h)+2\bbbz,$ then by \cite[Pro. 3.55(ii)]{Maz}, $W_1(n)$ is an irreducible $\fg$-module for which $e$ acts injectively. This is a contradiction as we have already assumed $\a\in R^{ln}$ and $e\in\LL^\a.$ So $$\tau({k_n},n)=(\mu+1)^2\hbox{ for some }\mu\in\lam(h)+2\bbbz.$$  Pick  $1\leq i\leq k$ such that $\lam(h)+2m_n$ is a weight  for  $W_i(n)/W_{i-1}(n)$ and note that  it is in fact a highest weight for $W_i(n)/W_{i-1}(n)$ due to Claim~1.  Since $W_i(n)/W_{i-1}(n)$ is irreducible,  it is a highest weight module of highest weight $\lam(h)+2m_n.$
So using \cite[Pro. 3.55(iii),(iv)]{Maz}, we get $(\lam(h)+2m_n+1)^2= \tau(k_n,n).$

\medskip

\noindent{\bf Claim~4.} There are positive integers $n_1,n_2,\ldots $ such that $\lam(h)+m_{n_1}$ is not a negative integer,
$$s_{n_1}\leq m_{n_1}<s_{n_2}
\leq m_{n_2}<s_{n_2}\leq \cdots$$ and
\begin{equation}\label{nonzero}
f^r (M^{\lam+m_{n_i}\a+(s_{n_i}+1)\b}\cap W^{(\tau(k_{n_i},n_i))})\neq \{0\}\quad\quad(r\in\bbbz^{>0}).
\end{equation}

\noindent{\it Reason:}
Recall (\ref{si}) and pick   $n_1\in\bbbz^{>0} $  in the way that   if $\lam(h)\in\bbbz,$ then  $-\lam(h)<s_{n_1}.$ So we get
\[\lam(h)+m_{n_1}\not\in \bbbz^{<0}.\]
Since $f$ acts   injectively, using (\ref{new}), we have
\[f^r( M^{\lam+m_{n_1}\a+(s_{n_1}+1)\b}\cap W^{(\tau(k_{n_1},n_1))})\neq \{0\}\quad\quad(r\in\bbbz^{>0}).\] We next  pick $n_2$ such that $s_{n_2}>m_{n_1}\geq s_{n_1}.$ As above, we have  \[f^r (M^{\lam+m_{n_2}\a+(s_{n_2}+1)\b}\cap W^{(\tau(k_{n_2},n_2))})\neq \{0\}\quad\quad(r\in \bbbz^{>0}).\]
Continuing  this process, we get the result.

\medskip

\noindent{\bf Claim~5.} $ \tau(k_{n_j},n_j)$'s are distinct.

\smallskip

\noindent{\it Reason:} Suppose to the contrary that for distinct $j,j',$ $\tau(k_{n_j},n_j)=\tau(k_{n_{j'}},n_{j'}).$ Then by Claim~3,  we have
\begin{align*}
\lam(h)^2+4m_{n_j}^2+1+&4m_{n_j}\lam(h)+4m_{n_j}+2\lam(h)=(\lam(h)+2m_{n_j}+1)^2\\=& \tau(k_{n_j},{n_j})\\
= &\tau(k_{n_{j'}},{n_{j'}})=(\lam(h)+2m_{n_{j'}}+1)^2\\=&\lam(h)^2+4m_{n_{j'}}^2+1+4m_{n_{j'}}\lam(h)+4m_{n_{j'}}+2\lam(h).
\end{align*}
So $4(m_{n_j}-m_{n_{j'}})(m_{n_j}+m_{n_{j'}}+\lam(h)+1)=0;$ that is $m_{n_j}+m_{n_{j'}}+\lam(h)+1=0.$ But this implies  that $\lam(h)\in\bbbz$ and  $$\lam(h)+m_{n_1}\leq \lam(h)+m_{n_j}=-(m_{n_{j'}}+1)\in \bbbz^{<0},$$ a contradiction; see Claim~4.

\medskip

\noindent{\bf Claim~6.} Recall $\b_1,\ldots,\b_t$ from the statement and for $1\leq i\leq t,$ fix $0\neq x_i\in \LL^{\b_i}.$ Consider the Casimir operator $c$ of $\fg$ as a linear transformation on $M$ and denote by $W^j$ ($j\in \bbbz^{>0}$), the generalized eigenspace of $c:M\longrightarrow M$  corresponding to $\tau(k_{n_j},n_j).$
i.e., $$W^j=\{v\in M\mid \exists r\in\bbbz^{>0}\ni \;(c-\tau(k_{n_j},n_j){\rm id})^rv=0\},$$ then $x_i^{s}v\in W^j$  for $s,j\geq 1,$  $1\leq i\leq t$ and $ v\in W^j.$

\smallskip

\noindent{\it Reason:}
  Since by our assumption, $(\a,\b_i)=0$ and  $\b_i\pm\a\not\in R,$ $c$  commutes with the action of $x_i$ on $M$; in particular,
  if $v\in M$ and $(c-\tau(k_{n_j},n_j){\rm id})^rv=0$ for some $r,j,$ then for $s\geq 1$ and   $1\leq i\leq t, $ we have
  \[(c-\tau(k_{n_j},n_j){\rm id})^rx_i^{s}v=x_i^{s}(c-\tau(k_{n_j},n_j){\rm id})^rv=0\] as we desired.

\medskip

\noindent{\bf Claim~7.}
Suppose $p\in\bbbz^{>1},$ then  $\dim(M^{\lam+s_{n_1}\a+(s_{n_{p}}+1)\b})\geq p;$ in particular, dimensions of weight spaces are not bounded:

\smallskip

\noindent{\it Reason:}
By our assumption, $-\a,\b_1,\ldots,\b_t\in R^{in}.$ So $f\in\LL^{-\a}$ as well as $0\neq x_1\in\LL^{\b_1},\ldots,0\neq  x_t\in\LL^{\b_t}$ act injectively on $M.$ Therefore, for $1\leq j\leq p,$ we have $s_{n_p}-s_{n_j},m_{n_{j}}-s_{n_{1}}\geq 0$  and contemplating (\ref{nonzero}), we have
\[ x_1^{s_{n_p}-s_{n_j}}\cdots x_t^{s_{n_p}-s_{n_j}} f^{m_{n_{j}}-s_{n_{1}}}(M^{\lam+m_{n_j}\a+(s_{n_j}+1)\b}\cap W^{(\tau(k_{n_j},n_j))})\neq \{0\}.\]
On the other hand, $W^{(\tau(k_{n_j},n_j))}$ is a $\fg$-module, so it is invariant under the action of $f.$ This together with Claim~6  and the fact  that  $W^{(\tau(k_{n_j},n_j))}\sub W^j,$
implies that
\[ x_1^{s_{n_p}-s_{n_j}}\cdots x_t^{s_{n_p}-s_{n_j}} f^{m_{n_{j}}-s_{n_{1}}} W^{(\tau(k_{n_j},n_j))}\sub W^{j}.\] These altogether imply that
\begin{align*}
&0\neq x_1^{s_{n_p}-s_{n_j}}\cdots x_t^{s_{n_p}-s_{n_j}} f^{m_{n_{j}}-s_{n_{1}}}(M^{\lam+m_{n_j}\a+(s_{n_j}+1)\b}\cap W^{(\tau(k_{n_j},n_j))})\\
&\sub  M^{\lam+s_{n_{1}}\a+(s_{n_{p}}+1)\b}\cap W^{j}.
\end{align*}
In other words, using Claim~5, the spaces $M^{\lam+s_{n_{1}}\a+(s_{n_{p}}+1)\b}\cap W^{j}$ $(1\leq j\leq p)$  are nonzero  linearly independent subspaces of $M^{\lam+s_{n_{1}}\a+(s_{n_{p}}+1)\b}$  which in turn implies that $\dim(M^{\lam+s_{n_{1}}\a+(s_{n_{p}}+1)\b})\geq p$ as we expected.
\qed
\begin{Pro}\label{final}
  Let $0\neq \dot\a,\dot\b\in\hbox{span}_\bbbr \dot R_{re}$  and $(\dot\a,\dot\b)=0$. Suppose there are  $\lam\in \supp(M),$ positive integers $n_1<n_2<n_3<\ldots$ and $t_1,t_2,\ldots\in \bbbz$ such that $\lam+4n_i(\dot\a+\dot\b)+t_i\d\in \supp(M)$ for all $i\geq 1.$ We have the following:
\begin{itemize}
\item[(i)] If $r\in\bbbz^{>0}$ and $k\in\{1,2\}$ with $k(-\dot\a+r \bbbz\d)\cup (2\dot\b+r \d+2r \bbbz\d )\sub R_{re},$  then
\begin{itemize}
\item[(a)] we have  either $k(-\dot\a+r \bbbz\d)\cap R^{ln}\neq \emptyset$ or
 $2\dot\b+r \d+2r \bbbz\d\sub R^{in},$

\item[(b)] we have either $(-2\dot\b-r \d+2r \bbbz\d)\cap  R^{ln}\neq \emptyset$ or $k(\dot\a+r \bbbz\d)\sub R^{in}.$
\end{itemize}

    \medskip

\item[(ii)] If $r\in\bbbz^{>0}$ and    $k\in\{1,2\}$  with $k\dot\a\pm 2\dot\b\not\in \dot R,$ then  the dimensions of weight spaces are not bounded provided that   the conditions of (a) or (b), as listed below, are satisfied:
\begin{itemize}
\item[(a)]
\begin{itemize}
\item[$\bullet$] $k(\dot\a+r \bbbz\d)\sub R^{ln}\cap R_0,$
\item[$\bullet$]     $k(-\dot\a+r \bbbz\d)\sub  R^{in}\cap R_0$ and
\item[$\bullet$]   there is $n^*\in\bbbz$ such that $2\dot\b+r \d+2r n^*\d\in R^{in},$
\end{itemize}
\medskip

\item[(b)]
\begin{itemize}
\item[$\bullet$]  $2\dot\b+r \d+2r \bbbz\d\sub R_0\cap R^{ln},$
\item[$\bullet$]      $-2\dot\b-r \d+2r \bbbz\d\sub R_0\cap R^{in}$ and
\item[$\bullet$]  there is $n^*\in\bbbz$ such that $k(\dot\a+r n^*\d)\in R^{in}.$
\end{itemize}
\end{itemize}
\end{itemize}
\end{Pro}
\pf Set
\[\dot\ep:=\dot\a+\dot\b\] and suppose $\lam\in \supp(M),$  $n_1<n_2<n_3<\cdots$ are positive integers  and $t_1,t_2,\ldots\in \bbbz$ such that
\begin{equation}\label{belong1}
\lam+4n_i\dot\ep+t_i\d\in \supp(M)\quad\quad (i\geq 1).
\end{equation}

For each $i\geq1,$ using division algorithm, we pick $q_i\in\bbbz$ and $d_i\in\{0,\ldots,2r -1\}$ such that  $t_i=2r q_i+d_i.$ Since $d_i$'s run over a finite set, infinitely many of them are equal; pick $i_1<i_2<\cdots$ such that  $d:=d_{i_1}=d_{i_2}=\cdots.$
Then
\begin{equation}\label{belong2}
\lam+4n_{i_j}\dot\ep+(2r q_{i_j}+d)\d=\lam+4n_{i_j}\dot\ep+t_{i_j}\d\stackrel{(\ref{belong1})}{\in} \supp(M)\quad\quad (j\geq 1);
\end{equation}
in particular,
\[
\mu:=\lam+4n_{i_1}\dot\ep+(2r q_{i_1}+d)\d\in \supp(M).
\]
Setting
$$m_j:=n_{i_{j+1}}-n_{i_1}\in\bbbz^{>0}\andd k_j:=q_{i_{j+1}}-q_{i_1}\quad(j\geq 1), $$ we have
\begin{align}
\mu+4m_j\dot\ep+2r k_j\d= & \mu+4(n_{i_{j+1}}-n_{i_1})\dot\ep+2r (q_{i_{j+1}}-q_{i_1})\d\nonumber \\
=&\lam+4n_{i_{j+1}}\dot\ep+(2r q_{i_{j+1}}+d)\d\stackrel{(\ref{belong2})}{\in}\supp(M)\quad\quad(j\geq 1).\label{in-supp2}
\end{align}

Next suppose $n^*\in\bbbz$ and    set
\begin{equation}\label{pj}
p_j:=k_j+m_j(-1-2n^*)\quad\quad (j\geq1).
\end{equation}
 Then we get for $j\geq1$ that
\begin{align}\label{first}
\hbox{\small$ \mu+4m_j\dot\a+2m_j(2\dot\b+r \d+2r n^*\d)+2r p_j\d$}
=&\hbox{\small $\mu+4m_j\dot\ep+2r (m_j+2n^*m_j+p_j)\d$}\nonumber\\
\stackrel{ (\ref{pj})}{=\joinrel=\joinrel=}&\hbox{\small $\mu+4m_j\dot\ep+2r k_j\d$}\stackrel{(\ref{in-supp2})}{\in} \supp(M)
\end{align}
and
\begin{align}\label{second}
\hbox{\small$ \mu+2m_j(2\dot\a+2r n^*\d)+2m_j(2\dot\b+r \d)+2r p_j\d
$}=&\hbox{\small $\mu+4m_j\dot\ep+2r (2n^*m_j+m_j+p_j)\d$}\nonumber\\
\stackrel{ (\ref{pj})}{=\joinrel=\joinrel=}&\hbox{\small $\mu+4m_j\dot\ep+2r k_j\d$}\stackrel{ (\ref{in-supp2})}{\in} \supp(M).
\end{align}

Now we are ready to go through the proof. We recall that $M$ has shadow and so
\[R_{re}^\times=R^{ln}\cup R^{in},\;\; R^{ln}\sub \frak{B}_M\andd R^{in}\sub \frak{C}_M.\]

(i)(a) To the contrary assume
\begin{equation}
\label{cont}
-k\dot\a+kr\bbbz\d\sub R^{in}\sub \frak{C}_M\andd \exists~ n^*\in\bbbz~~\ni 2\dot\b+r \d+2r n^*\d\in R^{ln}\sub \frak{B}_M.
\end{equation} Since $-k\dot\a+kr\bbbz\d\sub \frak{C}_M,$ recalling (\ref{pj}), we have $-k\dot\a-krp_j\d,-k{\dot\a}\in\frak{C}_M.$ But $k=1,2$ and $\frak{C}_M+\frak{C}_M\sub\frak{C}_M,$ so we get  $$-(2\dot\a+2r p_j\d),-2\dot\a\in \frak{C}_M\quad (j\geq 1).$$ Therefore, using (\ref{first}), we have
\begin{align*}
\mu+2m_j(2\dot\b+r \d+2r n^*\d)=&
\overbrace{\mu+2m_j(2\dot\b+r \d+2r n^*\d)+4m_j\dot\a+2r p_j\d}^{\in \supp(M)}\\-&(2m_j-1)2\dot\a-(2\dot\a+2r p_j\d)\in \supp(M)\quad\quad(j\geq 1)
\end{align*}
 which is a contradiction as by (\ref{cont}), $2\dot\b+r \d+2r n^*\d\in \frak{B}_M.$

\medskip
(i)(b) Suppose to the contrary that  $$-2\dot\b-r \d+2r \bbbz\d\sub  R^{in}\sub\frak{C}_M\andd \exists~ n^*\in\bbbz~~\ni k\dot\a+krn^*\d\in R^{ln}\sub\frak{B}_M.$$ Since $k=1,2,$  by \cite[Lem.~3.5]{you8},
\begin{equation}\label{cont2}
-2\dot\b-r \d+2r \bbbz\d\sub  \frak{C}_M\andd \exists~ n^*\in\bbbz~~\ni 2\dot\a+2r n^*\d\in \frak{B}_M.
\end{equation}

Therefore, recalling (\ref{pj}), we have $$-(2\dot\b+r \d+2r p_j\d),-(2\dot\b+r \d)\in \frak{C}_M\quad(j\geq 1).$$ So (\ref{second}) implies that
\begin{align*}
 \mu+2m_j(2\dot\a+2r n^*\d)=&
\overbrace{ \mu+2m_j(2\dot\a+2r n^*\d)+2m_j(2\dot\b+r \d)+2r p_j\d}^{\in\supp(M)}\\
 -&(2m_j-1)(2\dot\b+r \d)-(2\dot\b+r \d+2r p_j)\d\\
 &\in \supp(M)\quad\quad(j\geq 1)
 \end{align*} which is a contradiction as $2\dot\a+2r n^*\d\in \frak{B}_M;$ see (\ref{cont2}).
\medskip

(ii)(a)
Assume $r\in\bbbz^{>0},$ $k\in\{1,2\}$ and
\begin{equation}\label{con}
\hbox{\small $\begin{array}{ll}k\dot\a\pm 2\dot\b\not\in \dot R,&k(-\dot\a+r \bbbz\d )\sub  R^{in}\cap R_0\sub\frak{C}_M,\\
k(\dot\a+r \bbbz\d) \sub R^{ln}\cap R_0,& 2\dot\b+r \d+2r n^*\d\in R^{in}\sub\frak{C}_M,\end{array}$}
\end{equation} for some $n^*\in\bbbz.$ So recalling (\ref{pj}), we have  $$-k\dot\a,-(k\dot\a+krp_j\d)\in\frak{C}_M\quad\quad(j\geq1).$$  Therefore, by (\ref{first}), we have
\begin{align*}
 \hbox{\small $\mu+(2m_j-1)k\dot\a+2m_j(2\dot\b+r \d+2r n^*\d)$}=&
\hbox{\small$ \overbrace{\mu+4m_j\dot\a+2m_j(2\dot\b+r \d+2r n^*\d)+2r p_j\d}^{\in\supp(M)}$}\\
 -&\hbox{\small$(4/k-2)m_jk\dot\a-(k\dot\a+2r p_j\d)$}\\
 &\hbox{\small$\in  \supp(M)\quad\quad(j\geq 1).$}
\end{align*}
Setting $\a:=k\dot\a$ and $\b:=2\dot\b+r \d+2r n^*\d,$ we get using this together with  (\ref{con}) and the fact that $(\dot\a,\dot\b)=0$ that
\[\a\in R_0\cap R^{ln},\;-\a,\b\in R^{in},\;(\a,\b)=0,\;\a\pm\b\not\in R\] and
\[\mu+(2m_j-1)\a+2m_j\b\in \supp(M)\quad\quad(j\geq 1).\]
So Theorem \ref{BM} gives that dimensions of weight spaces are not bounded.

(ii)(b)
Suppose $r\in\bbbz,$  $k\in\{1,2\},$
\begin{equation}\label{con2}
\hbox{\small$
\begin{array}{ll}
2\dot\b+r \d+2r \bbbz\d\sub R_0\cap R^{ln},&
  -2\dot\b-r \d+2r \bbbz\d\sub R_0\cap R^{in}\sub\frak{C}_M\\
2\dot\b\pm k\dot\a\not\in\dot R,& k\dot\a+krn^*\d\in R^{in}\sub\frak{C}_M,
 \end{array}$}
 \end{equation}
  for some  $n^*\in\bbbz.$
For $p_j$ as in (\ref{pj}), since $-2\dot\b-r \d-2r p_j\d\in \frak{C}_M,$ (\ref{second}) implies that
\begin{align*}
&\hbox{\small $\mu+(2m_j-1)(2\dot\b+r \d)+2m_j(2\dot\a+2r n^*\d)$}\\
=&\hbox{\small $\underbrace{\mu+2m_j(2\dot\a+2r n^*\d)+2m_j(2\dot\b+r \d)+2r p_j\d}_{\in\supp(M)}$}-\hbox{\small $(2\dot\b+r \d+2r p_j\d)$}\in \supp(M).
\end{align*} This together with (\ref{con2}) and the fact that $(\dot\a,\dot\b)=0,$ implies that for
\[\a:=2\dot\b+r \d\andd \b:=\left\{
\begin{array}{ll}
\dot\a+r n^*\d& k=1,\\
2\dot\a+2r n^*\d& k=2,
\end{array}
\right.
\] we have
\[\a\in R_0\cap R^{ln},\;-\a,\b\in R^{in},\;(\a,\b)=0,\;\a\pm\b\not\in R\] and
\[\left\{
\begin{array}{ll}
\mu+(2m_j-1)\a+2m_j\b\in \supp(M)&k=2,\\
\mu+(2m_j-1)\a+2m_j(\b+\b)\in \supp(M)&k=1,\\
\end{array}
\right.
\]
for $j\geq 1.$
So  Theorem~\ref{BM} gives what we desired.
\qed

\begin{cor}\label{final-2}
 For $\dot\gamma\in \dot R\setminus\{0\},$ recall $S_{\dot\gamma}$  from (\ref{s})  and assume $\boldsymbol\zeta:\hbox{\rm span}_\bbbr \dot R\longrightarrow \bbbr$ is a linear functional such that
\begin{equation}
\label{pos}
\boldsymbol\zeta(\dot\gamma)>0\Leftrightarrow\dot\gamma+S_{\dot\gamma}\sub R^{\ln}\andd -\dot\gamma+S_{-\dot\gamma}\sub R^{in}\quad\quad(\dot\gamma\in {\dot R_{re}}\setminus\{0\}).
\end{equation}
Suppose there are  $\dot\ep\in \dot R_{ns}\setminus\{0\}$ with   $\boldsymbol\zeta(\dot\ep)>0,$   $\lam\in\supp(M)$ and $k_1,k_2,\ldots\in \bbbz$ such that $$\lam+4n\dot\ep+k_n\d\in \supp(M)\quad\quad(n\in\bbbz^{>0}).$$ Then
dimensions of weight spaces are unbounded.
\end{cor}
 \pf
Recalling $r_{\dot\ep}$ from (\ref{salpha}) and using Remark~\ref{rem}, we get $k=1,2$ and  $\dot \a,\dot\b$ such that
\begin{equation}\label{sub}k(\pm\dot\a+r_{\dot\ep}\bbbz\d)\cup (\pm2\dot\b+r_{\dot\ep}\d+2r_{\dot\ep}\bbbz\d )\sub R_0\cap R_{re} \andd \dot\ep=\dot\a+\dot\b
\end{equation} and $k\dot\a\pm2\dot \b\not\in \dot R.$
Since $\boldsymbol\zeta(\dot\a)+\boldsymbol\zeta(\dot\b)=\boldsymbol\zeta(\dot\ep)>0,$  the value of $\boldsymbol\zeta$ on at least one of $\dot\a$ and $\dot\b$ is positive.

\noindent{\bf Case~1.} $\boldsymbol\zeta(\dot\a)>0,$ then (\ref{pos}) and (\ref{sub}) imply that  $k(\dot\a+r_{\dot\ep}\bbbz\d)\sub R^{ln}\cap R_0$ and
$k(-\dot\a+r_{\dot\ep}\bbbz\d)\sub R^{in}\cap R_0.$  Using Proposition~\ref{final}(i)(a), we get $2\dot\b+r_{\dot\ep}\d+2r_{\dot\ep}\bbbz\d \sub R^{in}.$ So the conditions of Proposition~\ref{final}(ii)(a) are fulfilled and so we are done.

\noindent{\bf Case~2.} $\boldsymbol\zeta(\dot\b)>0,$ then by (\ref{pos}) and (\ref{sub}), we have  $2\dot\b+r_{\dot\ep}\d+2r_{\dot\ep}\bbbz\d \sub R^{ln}\cap R_0$ and $-(2\dot\b+r_{\dot\ep}\d+2r_{\dot\ep}\bbbz\d )\sub R^{in}.$ So by Proposition~\ref{final}(i)(b), we have $k(\dot\a+r_{\dot\ep}\bbbz\d)\sub R^{in}.$  Therefore, we get the result using  Proposition~\ref{final}(ii)(b).
\qed
\section{Tight irreducible   finite weight modules}
Irreducible finite weight  modules over affine Lie superalgebras are divided into two classes hybrid and tight. In \cite{you8}, we have studied hybrid modules and here we start studying tight modules over twisted affine Lie superalgebras with nonzero odd part. Throughout this section, we assume $\LL$ is one of the twisted affine Lie superalgebras $A(2k-1,2\ell-1)^{(2)}$ ({\tiny $(k,\ell)\neq(1,1)$}), $A(2k,2\ell)^{(4)},$ $A(2k,2\ell-1)^{(2)}$ and $D(k+1,\ell)^{(2)}$ where $k$ is a nonnegative integer and $\ell$ is a positive integer.

\begin{Thm}[{\cite[Thm.~4.8]{you8}}]\label{property}
Suppose that $M$ is a {weight}  $\LL$-module having shadow, then for each $\b\in R_{re}^\times,$ one of the following will happen:
\begin{itemize}
\item[\rm (i)] $(\b+\bbbz\d)\cap R\sub R^{ln},$
\item[\rm (ii)] $(\b+\bbbz\d)\cap R\sub R^{in},$
\item[\rm (iii)] there exist $m\in\bbbz$ and $t\in\{0,1,-1\}$ such that for $\gamma:=\b+m\d,$
\begin{align*}
&(\gamma+\bbbz^{\geq 1}\d)\cap R\sub R^{in},\quad (\gamma+\bbbz^{\leq 0}\d)\cap R\sub R^{ln}\\
& (-\gamma+\bbbz^{\geq t}\d)\cap R  \sub R^{in},\quad (-\gamma+\bbbz^{\leq t-1}\d)\cap R\sub R^{ln},
\end{align*}
\item[\rm (iv)]there exist $m\in\bbbz$ and $t\in\{0,1,-1\}$ such that for $\eta:=\b+m\d,$
\begin{align*}
&(\eta+\bbbz^{\leq -1}\d)\cap R\sub R^{in},\quad (\eta+\bbbz^{\geq 0}\d)\cap R\sub R^{ln}\\
& (-\eta+\bbbz^{\leq -t}\d)\cap R  \sub R^{in},\quad (-\eta+\bbbz^{\geq 1-t}\d)\cap R\sub R^{ln}.
\end{align*}
\end{itemize}
\end{Thm}

\begin{deft}\label{full-hybrid}
{\rm Suppose that $M$ is a weight $\LL$-module having shadow. {We say that} $\a\in R_{re}^\times$ is {\it full-locally nilpotent} (resp. {\it full-injective}) if $(\a+\bbbz\d)\cap R\sub R^{ln}$  (resp. $(\a+\bbbz\d)\cap R\sub R^{in}$), otherwise, we call it {\it hybrid}. The module $M$ is called {\it hybrid} if all nonzero real  roots are hybrid and otherwise it is called {\it tight}.
}
\end{deft}
For a weight  $\LL$-module $M$  having shadow, denote the set of all  full-locally nilpotent (resp. full-injective) real roots by $R_{\rm f-ln}$ (resp. $R_{\rm f-in}$) and the set of all hybrid real roots by $R_{\rm hyb}.$
\begin{rem}
\label{rem2}
{\rm
Suppose that $M$ is a weight $\LL$-module having shadow and $\a$ is a nonzero real odd root such that $2\a\in R_{\rm hyb}.$ We claim that $\a\in R_{\rm hyb}$.
If $\a\not\in R_{\rm hyb},$ then either $\a\in R_{\rm f-ln}$ or $\a\in R_{\rm f-in}.$ Suppose $r$ is as in (\ref{r}), then  either $\a+r\bbbz\d \sub R^{ln}$ or $\a+r\bbbz\d\sub  R^{in},$ respectively. Contemplating (\ref{equal}), in the former case, we have $2\a+2r\bbbz\d \sub R^{ln}$ and in the latter case, we have  $2\a+2r\bbbz\d\sub  R^{in}$ due to \cite[Lem.~3.5]{you8} and the fact that $M$ has shadow but both give contradictions as $2\a\in R_{\rm hyb}$. }
\end{rem}

\begin{lem}\label{fini}
Suppose  $\D$ is an irreducible finite root system with the inner product  $\fm$ on its real linear span. Denote the length of a root $\a$ by $\ell(\a).$ If $\a,\b\in \D\setminus\{0\}$ with $\a+\b\in \D\setminus\{0\},$ then one of the following occurs:
\begin{itemize}
\item[(a)] $\ell(\a)=\ell(\b)<\ell(\a+\b),$
\item[(b)] $\ell(\a+\b)=\ell(\a)<\ell(\b)$ or  $\ell(\a+\b)=\ell(\b)<\ell(\a),$
\item[(c)] $\ell(\a)=\ell(\b)=\ell(\a+\b).$
\end{itemize}
\end{lem}
\pf {It is easily verified.}
\qed

\medskip

The following proposition is the super version of  \cite[Pro.~2.22-Case~2]{DG}.
\begin{Pro}\label{para}
 Recall $R_0(i)$ and $\dot R_0(i)$ $(i=1,2)$ from Tables~\ref{table3},\ref{table4} and let  $M$ be a weight  $\LL$-module having shadow.  Suppose that
\begin{itemize}
\item   $R_{re}^\times\neq R_{hyb}$,
\item $R^{ln}\cap R_0(1)$  and $R^{ln}\cap R_0(2)$ (if $R_0(2)$ is nonempty) are  nonempty proper subsets of $R_0(1)\cap R_{re}$ and $R_0(2)\cap R_{re}$ respectively.
\end{itemize}
Set
$$P:=R_{\rm f-ln}\cup-R_{\rm f-in}\cup R_{\rm hyb}\cup \bbbz\d.
$$
 Suppose $i=1,2$  and  $R_0(i)\neq\emptyset,$ then we have the following:
\begin{itemize}
\item[(i)] $P$ is a closed subset of $R_{re}\cup R_{im}$ with $P\cup-P=R_{re}\cup R_{im};$ in particular,  $P_i:=P\cap R_0(i)$ is a parabolic subset of $R_0(i).$
\item[(ii)] $\dot P_i:=\{\dot\a\in \dot R_0(i)\mid \exists m\in\bbbz\; \hbox{ \rm  s.t. } \dot\a+m\d\in P_i\}$    is a parabolic subset  of $\dot R_0(i)$ satisfying  $P_i=P\cap R_0(i)=(\dot P_i+\bbbz\d)\cap R_0(i).$ Moreover, at least one of $\dot P_1$ and $\dot P_2$ is proper.
\end{itemize}
\end{Pro}
\pf
(i) {It  follows easily} from \cite[Lem.~3.5 \& Thm.~4.7]{you8} and Theorem~\ref{property} but for the convenience of readers, we prove one case {to explain the argument.} Suppose that $\a\in R_{f-ln}$ and $\b\in R_{hyb}$ (equivalently $-\b\in R_{hyb}$ by Theorem~\ref{property}) with $\gamma:=\a+\b\in R_{re}.$ We claim that $\a+\b\in R_{f-ln}.$ Suppose  to the contrary, $\a+\b\not\in R_{f-ln}.$
Let $r$ be as in (\ref{r}). Then using (\ref{equal}), we have \[\gamma+r\bbbz\d\sub R_{re}^\times\quad\quad (\gamma\in R_{re}^\times).\]
Since $\a+\b\not\in R_{f-ln},$ then $\a+\b+m\d\in R_{in}$  for some $m\in r\bbbz.$ Since $-\b\in R_{hyb},$ we pick an integer $n\in r\bbbz$ such that $-\b+n\d\in R_{in}.$ Then, since $\frak{C}_M+\frak{C}_M\sub\frak{C}_M$ and $R_{re}^\times\cap \frak{C}_M=R_{in},$ we have $$\underbrace{\a+(m+n)\d}_{R_{re}^\times}=(\underbrace{\a+\b+m\d}_{\in R_{in}\sub\frak{C}_M})+(\underbrace{-\b+n\d}_{\in R_{in}\sub\frak{C}_M})\in R_{re}^\times\cap \frak{C}_M= R_{in}$$ which is a contradiction.

 (ii)
We suppose $\dot\a_1,\dot\a_2\in \dot P_i$ and $\dot\a_1+\dot\a_2\in \dot R_0(i).$ We need to show that $\dot\a_1+\dot\a_2\in \dot P_i.$ We recall subsets $$S_{\dot\a_j}(i)=\{m\d\in\bbbz\d\mid \dot\a_j+m\d\in R_0(i)\}\quad\quad(j=1,2)$$  from  (\ref{s-i}). Since $\dot\a_1,\dot\a_2\in \dot P_i,$ there are $\sg_1\in S_{\dot\a_1}(i)$ and $\sg_2\in S_{\dot\a_2}(i)$ with $$\dot\a_1+\sg_1\in P_i=R_0(i)\cap P\andd \dot\a_2+\sg_2\in P_i=R_0(i)\cap P.$$
Since $\pm\d\in P,$ if $j=1,2$ and  $\sg\in S_{\dot\a_j}(i),$ we have
\[\underbrace{\dot\a_j+\sg}_{R_{re}\cup R_{im}}=\underbrace{\dot\a_j+\sg_j}_{P}+\underbrace{\sg-\sg_j}_{P},\]
but $P$ is a closed subset of  $R_{re}\cup R_{im},$ so $\dot\a_j+\sg\in P.$ Therefore,  we get
\begin{equation}\label{sub P}
\dot\a_1+S_{\dot\a_1}(i)\sub  P_i=R_0(i)\cap P\andd \dot\a_2+S_{\dot\a_2}(i)\sub  P_i=R_0(i)\cap P.
\end{equation}

Without  loss of generality, we assume $\dot\a_1,\dot\a_2,\dot\a_3:=\dot\a_1+\dot\a_2\neq 0.$ By our assumption,  $\dot\a_1+\dot\a_2\in  \dot R_0(i).$
One  knows from Lemma~\ref{fini} that there are ${s'}\in\{1,2\},$ $s,k\in\{1,2,3\}$ with $\{1,2,3\}=\{s,{s'},k\}$ and $t=\left\{
\begin{array}{ll}
1&s\neq 3\\
-1&s=3
\end{array}
\right.
$  such that
\[\dot\a_k=\dot\a_s+t\dot\a_{s'}\andd \ell(\dot\a_s)=\ell(t\dot\a_{s'})\leq \ell(\dot\a_k).\]
So by (\ref{sum}), we have
\begin{equation}\label{subsub}
S_{\dot\a_k}(i)\sub S_{\dot\a_s}(i)+S_{t\dot\a_{s'}}(i).
\end{equation}

If $s\neq 3,$ then $k=3$ and  $t=1$. For $\sg\in S_{\dot\a_3}(i)=S_{\dot\a_k}(i),$ using (\ref{subsub}), one finds $\tau\in S_{\dot\a_s}(i)$ and $\tau'\in S_{\dot\a_{s'}}(i)=S_{t\dot\a_{s'}}(i)  $ with $\sg=\tau+\tau'$ and so by (\ref{sub P}), we have
\[(\dot\a_1+\dot\a_2)+\sg=\dot\a_3+\sg=(\dot\a_s+\tau)+(\dot\a_{s'}+\tau')\in R_0(i)\cap ( P_i+ P_i)\sub P_i \]
which in turn implies that $$\dot\a_1+\dot\a_2\in \dot P_i.$$

Also if  $s=3,$ then $t=-1.$ Pick $\tau\in S_{\dot\a_k}(i).$ Then, by (\ref{subsub}), there are $\sg\in S_{\dot\a_s}(i)$ and   $\gamma\in S_{-\dot\a_{s'}}(i)$  (equivalently, $-\gamma\in S_{\dot\a_{s'}}(i)$) with $\tau=\sg+\gamma,$ so using (\ref{sub P}), we get
\[(\dot\a_1+\dot\a_2)+\sg=\dot\a_s+\sg=(\dot\a_k+\tau)+(\dot\a_{s'}-\gamma)\in R_0(i)\cap ( P+ P)\sub P. \]
which in turn implies that $$\dot\a_1+\dot\a_2\in \dot P_i.$$
So $\dot P_i$ is a parabolic subset of $\dot R_0(i)$. {Finally, assume that neither  $\dot P_1$ nor $\dot P_2$
is  proper}. {As for $i=1,2,$ $R_0(i)\cap R^{ln}$ is a nonempty proper subset of $R_0(i)\cap  R_{re}$,} the  same argument as in \cite[Pro.~2.22-Case~2]{DG}, gives that for $i=1,2,$ $R_0(i)\cap R_{re}^\times\sub R_{hyb}.$
This together with Remark~\ref{rem2} implies that $R_{re}^\times=R_{hyb},$ {which is a contradiction}. So at least  one of $\dot P_1$ and $\dot P_2$ is proper.
\qed

\medskip

\begin{Thm}\label{main-2}
Suppose that $M$ is a tight irreducible finite weight $\LL$-module with bounded weight multiplicities such that
\begin{itemize}
\item $R^{ln}\cap R_0(1)$ is a  nonempty proper subset of  $R_0(1)\cap R_{re}$ and if $R_0(2)$ is nonempty, $R^{ln}\cap R_0(2)$ is a nonempty proper subset of   $R_0(2)\cap R_{re}$.
\end{itemize}
Then there is  a nontrivial  triangular decomposition  $$R=R^+\cup R^\circ\cup R^-$$ for $R$ such that
$M^{\LL^+}=\{v\in M\mid \LL^\a v=\{0\}\;\;(\forall \a\in R^+)\}\neq\{0\}.$ In particular, $M^{\LL^+}$ is an irreducible finite weight module over $\LL^\circ=\sum_{\a\in R^\circ}\LL^\a$ and $M\simeq {\rm Ind}_\LL(M^{\LL^+}).$
\end{Thm}
\pf
Recall $P,$ $P_i$ and $\dot P_i$ ($i=1,2$) from Proposition~\ref{para}. If $R_0(i)$ ($i=1,2$) is nonempty, then $\dot P_i$ is a parabolic subset of $\dot R_0(i),$ so by \cite[Pro. 2.10]{DFG}, there is
$\boldsymbol{\zeta_i}:\sspan_\bbbr \dot R_0(i)\longrightarrow \bbbr$ such that
$$\dot P_i=\{\dot\a\in \dot R_0(i)\mid {\boldsymbol\zeta_i}(\dot\a)\geq 0\}.$$
Define $${\boldsymbol\zeta}:=\left\{
\begin{array}{ll}
{\boldsymbol\zeta_1}\op\boldsymbol{\zeta_2}:\sspan_\bbbr \dot R_0(1)\op \sspan_\bbbr \dot R_0(1)\longrightarrow \bbbr& R_0(2)\neq \emptyset\\
\boldsymbol{\zeta_1}:\sspan_\bbbr  \dot R_0(1)\longrightarrow \bbbr& R_0(2)= \emptyset.
\end{array}
\right.
$$
Extend $\boldsymbol\zeta$ to $\sspan_\bbbr R_0=\sspan_\bbbr R$ (see (\ref{span})) with $\boldsymbol{\zeta}(\d)=0.$

Recall (\ref{s}) and  suppose $\dot\a\in \dot R_{re}^\times$ and $m\d\in S_{\dot\a}.$ Since two times of a real odd root $\a$ (i.e., real root $\a$ with  $\a\in R_1=R\setminus R_0$) is a real even root, depending on $\a:=\dot\a+m\d\in R_0$ or $\a=\dot\a+m\d\in R_1,$  there are $i=1,2$ and $k=1$ or $k=2$ such that $\b:=k(\dot\a+m\d)\in R_0(i).$  Since by Proposition~\ref{para}, $P$ is a closed subset of $R_{re}\cup R_{im}$ with  $P\cup -P=R_{re}\cup R_{im},$  we have $$\a\in P\Leftrightarrow\b\in P_i\Leftrightarrow k\dot\a\in\dot P_i\Leftrightarrow \boldsymbol{\zeta}(\dot\a)\geq 0\Leftrightarrow \boldsymbol{\zeta}(\a)\geq 0.$$
In particular,
\[\boldsymbol{\zeta}(\dot\a)>0\Leftrightarrow \dot\a+S_{\dot\a}\sub R^{ln}\andd -\dot\a+S_{-\dot\a}\sub R^{in}\quad\quad(\dot \a\in \dot R_{re}\setminus\{0\}).\] We mention that by Proposition~\ref{para}(ii) either $\dot P_1$ is proper or $\dot P_2$ is proper, so $\boldsymbol\zeta$ is nonzero.

 For the triangular decomposition corresponding to $\boldsymbol{\zeta},$
 use Proposition \ref{pdelta} for the $\LL$-module $M$ to choose $\lam\in\supp(M)$ such that $\lam+\a\not\in \supp(M)$ for all $\a\in R_{re}^+.$ This  in particular implies that  $$\LL^\a M^\lam=\{0\}\quad\quad(\a\in R_{re}^+).$$ So
$$B:=\{v\in M\mid \hbox{ $v$ is a weight vector and }\LL^\a v=\{0\}\;\;(\forall \a\in R_{re}^+)\}$$ is a nonzero subspace of $M.$
For $v\in B,$ set
$$\aa_v:=\{\dot\ep\in \dot R_{ns}^\times\mid  \exists~ m\in \bbbz \hbox{\; s.t. \;}  \LL^{\dot\ep+m\d}v\neq \{0\}\hbox{ and $\boldsymbol{\zeta}(\dot\ep+m\d)=\boldsymbol{\zeta}(\dot\ep)>0$}\}.$$

\noindent{\bf Step 1.}
Suppose  $v\in B$ is such that $0\neq |A_v|={\rm min}\{|A_{x}|\mid x\in B\}$  and  $\dot\ep\in \aa_v$ is such that $\boldsymbol\zeta(\dot\ep)={\rm max}\{\boldsymbol\zeta(\dot\eta)\mid \dot\eta\in \aa_v \}.$ If $n\in \bbbz$ and $0\neq w\in \LL^{\dot \ep+n\d}v,$ then, we have $$w\in B,\;\;\;\aa_w=\aa_{v}\andd \boldsymbol\zeta(\dot\ep)={\rm max}\{\boldsymbol\zeta(\dot\eta)\mid \dot\eta\in \aa_w \}:$$
{ To show that $w\in B,$ we need to prove $\LL^\a w=\{0\}$ for all $\a\in R_{re}^+.$ So suppose $\a\in R_{re}^+.$
We have}
$$\LL^\a w\sub \LL^\a\LL^{\dot\ep+n\d}v\sub \LL^{\dot\ep+n\d}\underbrace{\LL^\a v}_0+[\LL^{\dot\ep+n\d},\LL^\a ]v\sub \LL^{\dot\ep+n\d+\a}v.$$  {So it is enough to show that $\LL^{\dot\ep+n\d+\a}v=\{0\}.$ To this end, we consider the following:}
\begin{itemize}
\item $\dot\ep+n\d+\a\not \in R:$  In this case, we get $\LL^\a w\sub  \LL^{\dot\ep+n\d+\a}v=\{0\}.$
\item $\dot\ep+n\d+\a \in R_{re}:$ Since     $\a\in R_{re}^+,$ we have $\boldsymbol\zeta(\a)>0$ and so $\boldsymbol\zeta(\dot\ep+n\d+\a)=\underbrace{\boldsymbol\zeta(\dot\ep)}_{>0}+
    \underbrace{\boldsymbol\zeta(\a)}_{>0}+n\underbrace{\boldsymbol{\zeta}(\d)}_{=0}
    >0.$ Therefore,     we get that $\LL^\a w\sub  \LL^{\dot\ep+n\d+\a}v\stackrel{v\in B}{=\joinrel=}\{0\}.$
\item   $\dot\ep+n\d+\a \in R_{ns}:$  Let $\a=\dot\a+\sg$ for some $\dot \a\in \dot R$ and $\sg\in R^0.$ Since $\dot\ep+n\d+\a \in R_{ns},$ we have $\dot\ep+\dot\a\in \dot R_{ns}$ and since $\a\in R_{re}^+$ and $\boldsymbol{\zeta}(\d)=0,$ we have $\boldsymbol{\zeta}(\dot\a)>0.$ But   $\boldsymbol\zeta(\dot\ep)={\rm max}\{\boldsymbol\zeta(\dot\eta)\mid \dot\eta\in \aa_v \}$ while   \[\boldsymbol\zeta(\dot\ep+\dot\a)=
    \boldsymbol\zeta(\dot\ep)+\underbrace{\boldsymbol\zeta(\dot\a)}_{>0}>\boldsymbol{\zeta}(\dot\ep),\] so  we get that $\dot\ep+\dot\a\not \in\aa_v$ and so $\LL^\a w\sub  \LL^{\dot\ep+n\d+\a}v=\{0\}.$
\end{itemize}
     These altogether imply that $w\in B.$
Next suppose  $\dot\eta\in \aa_{w}.$ Then $\boldsymbol\zeta(\dot\eta)>0$ and  there is $m\in\bbbz$ such that $\LL^{\dot\eta+m\d}w\neq \{0\}.$
So \begin{align*}
\{0\}\neq \LL^{\dot\eta+m\d}w\sub \LL^{\dot\eta+m\d}\LL^{\dot\ep+n\d}v\sub& [\LL^{\dot\eta+m\d},\LL^{\dot\ep+n\d}]v +\LL^{\dot\ep+n\d}\LL^{\dot\eta+m\d} v\\
\sub & \LL^{\dot\ep+\dot\eta+m\d+n\d}v +\LL^{\dot\ep+n\d}\LL^{\dot\eta+m\d} v.
\end{align*}
But if $\dot\ep+\dot\eta+m\d+n\d\not\in R,$ we have $\LL^{\dot\ep+\dot\eta+m\d+n\d}v=0.$ Also if $\dot\ep+\dot\eta+m\d+n\d\in R,$ since  by (\ref{re-re}), $(R_{ns}^\times+R_{ns}^\times)\cap R\sub R_{re}\cup R_{im}$ and  $\boldsymbol\zeta(\dot\ep+\dot\eta)=\boldsymbol\zeta(\dot\ep)+\boldsymbol\zeta(\dot\eta)>0,$ we get that
$\dot\ep+\dot\eta+m\d+n\d\in R_{re}^+$ and so again we get $\LL^{\dot\ep+\dot\eta+m\d+n\d}v=0.$ Therefore, we have
\begin{align*}
\{0\}\neq \LL^{\dot\eta+m\d}w\sub \LL^{\dot\ep+n\d}\LL^{\dot\eta+m\d} v.
\end{align*}
This implies that $\LL^{\dot\eta+m\d} v\neq \{0\}.$ Therefore, $\dot\eta\in \aa_{v},$ in other words $\aa_{w}\sub \aa_{v}.$ But $\aa_v$ has {the smallest cardinality}, so  we get   $\aa_{w}=\aa_{v}.$ This completes the proof of this step.

\noindent{\bf Step 2.}
Pick $v_0\in B$ such that $\aa_{v_0}$ is of the smallest cardinality, then $\aa_{v_0}=\emptyset:$ To the contrary, assume $\aa_{v_0}\not =\emptyset$ and  pick $\dot\ep\in \aa_{v_0}$ such that $\boldsymbol\zeta(\dot\ep)={\rm max}\{\boldsymbol\zeta(\dot\eta)\mid \dot\eta\in \aa_v \}.$ Since $v_0\in B,$  there are $\lam\in\supp(M)$ with $v_0\in M^\lam$ and  since $\dot \ep\in \aa_{v_0},$ there is  $m_1\in\bbbz$ with  $ \LL^{\dot\ep+m_1\d}v_0\neq\{0\}.$ Pick $0\neq v_1\in \LL^{\dot\ep+m_1\d}v_0,$ then by Step~1, we have
$$v_1\in B,\;\;\;\aa_{v_1}=\aa_{v_0}\andd \boldsymbol\zeta(\dot\ep)={\rm max}\{\boldsymbol\zeta(\dot\eta)\mid \dot\eta\in \aa_{v_1} \};$$ in particular $\aa_{v_1}$ is of {the smallest cardinality}. Since $v_1\in B$ and $\dot\ep\in \aa_{v_0}=\aa_{v_1},$ we pick $m_2\in\bbbz^{>0}$ and $0\neq v_2\in \LL^{\dot\ep+m_2\d}v_1,$ again using Step~1, we have
$$v_2\in B,\;\;\;\aa_{v_2}=\aa_{v_1}=\aa_{v_0}\andd \boldsymbol\zeta(\dot\ep)={\rm max}\{\boldsymbol\zeta(\dot\eta)\mid \dot\eta\in \aa_{v_2} \}.$$

Continuing  this process, one finds integers $m_1,m_2,m_3,\ldots$ and nonzero vectors $v_1,v_2,v_3,\ldots$ with $0\neq v_i\in\LL^{\dot\ep+m_i\d}v_{i-1}.$ In particular, $$\lam+n\dot\ep+(m_1+\cdots+m_n)\d\in\supp(M)\quad\quad(n\geq 1).$$ Setting   $k_{n}:=m_1+\cdots+m_{4n},$ we have  $$ \lam+4n\dot\ep+k_n\d\in\supp(M)\quad\quad(n\geq 1). $$ This together with Corollary~\ref{final-2} and the fact that  $\boldsymbol\zeta(\dot\ep)>0$ gives a contradiction. Therefore $\aa_{v_0}=\emptyset$ and so we are done using  Proposition \ref{ind}(ii). \qed

\smallskip

\centerline{\bf Acknowledgment.}
The author would like to thank the anonymous referee for his/her fruitful comments and suggestions.

\end{document}